\documentclass[10pt, twocolumns]{article}

\usepackage[a4paper,margin=1.7cm]{geometry}
\usepackage[T1]{fontenc}
\usepackage[utf8]{inputenc}
\usepackage{lmodern}

\usepackage{amsmath,amssymb,amsfonts,mathtools,bm}
\usepackage{graphicx}
\usepackage{booktabs}
\usepackage{multirow}
\usepackage{array}
\usepackage{xcolor}
\usepackage{hyperref}
\usepackage{enumitem}
\usepackage{setspace}
\usepackage{caption}
\usepackage{subcaption}
\usepackage{dblfloatfix}
\usepackage{float}
\usepackage[section]{placeins}
\usepackage{microtype}

\hypersetup{
    colorlinks=true,
    linkcolor=blue,
    citecolor=blue,
    urlcolor=blue
}

\title{\textbf{Multimodal branched transport infers anatomically aligned brain reaction maps}}

\author{
Cristian Mendico\\
\small Institut de Math\'ematique de Bourgogne, UMR 5584 CNRS, Universit\'e Bourgogne Europe\\
\small \texttt{cristian.mendico@u-bourgogne.fr}
}
\date{}

\begin{document}

\twocolumn[
\maketitle

\begin{abstract}
How external stimulation is transformed into distributed reaction patterns remains unresolved at the level of propagation architecture. Existing large-scale control models quantify transition costs on prescribed networks but do not infer the routing map itself from source and target activity. Here we combine task-related blood-oxygen-level-dependent responses, source-reconstructed electrophysiology and tractography-derived anisotropy to estimate stimulation and reaction measures, define an anatomical transport cost, and infer a branched propagation architecture by variational optimisation. Unlike standard transport formulations, branched transport favours aggregation of signal into shared neural highways before redistribution. We further attach a stochastic graph-induced dynamics to the inferred map and quantify the trade-off between geometric efficiency and dynamical controllability. We show that multimodal data generate anatomically aligned brain reaction maps, that anisotropic costs qualitatively reshape routing backbones relative to isotropic baselines, and that hybrid geometric--dynamical optimisation reveals non-trivial rank reversals across branching regimes.

\smallskip 

\end{abstract}

\vspace{0.5cm}
]

\section*{Introduction}

A central problem in systems neuroscience is to determine how an external stimulation is propagated through the brain so as to produce a reaction. Current deterministic and stochastic control approaches quantify the cost of driving brain states on a prescribed structural substrate, but they do not treat the routing architecture itself as an unknown \cite{Gu2015,Kawakita2022,Kamiya2023}. In the companion theoretical paper, we proposed a different variational viewpoint: the unknown is a graph or current connecting a stimulation source measure to a reaction target measure, selected as the minimiser of an anisotropic branched transport problem \cite{MendicoTheory2026}. In that framework, the support of the minimiser defines a \emph{brain reaction map}, namely an inferred stimulus-to-reaction routing architecture rather than a controlled trajectory on a fixed network \cite{MendicoTheory2026}.

The present article develops the data-driven computational counterpart of that theory. Instead of assuming abstract source and target measures, we estimate them from multimodal neuroimaging signals. Instead of postulating a generic transport metric, we derive a geometry-aware cost from tractography-informed anisotropy \cite{Hagmann2007,Yeh2021,Zhang2022}. Instead of stopping at the geometric minimiser, we use the inferred graph as the substrate of a stochastic whole-brain dynamics and evaluate the trade-off between anatomical transport efficiency and dynamical control cost. This places the work at the intersection of multimodal neuroimaging, connectomics and large-scale brain modelling \cite{Breakspear2017,Deco2020,Gilson2020,Patow2024,Pathak2022}.

This leads to three concrete questions. First, can task-related blood-oxygen-level-dependent responses and source-reconstructed EEG/MEG be fused into coherent source and target measures for a transport problem? Second, does tractography-derived anisotropy substantially alter the inferred routing architecture relative to isotropic baselines? Third, once a graph is inferred, is the geometrically optimal map also dynamically plausible as a substrate for stochastic state transitions?

To answer these questions, we implement the full multimodal pipeline already sketched in the theoretical work \cite{MendicoTheory2026}: (i) estimation of $\mu_{\mathrm{stim}}^{+}$ and $\mu_{\mathrm{react}}^{-}$ from fMRI and EEG/MEG, (ii) construction of an anisotropic transport cost from diffusion-informed tensors, (iii) branched transport optimisation on a candidate graph, and (iv) a hybrid stochastic extension built on the inferred routing architecture. The resulting article is intentionally complementary to the theoretical companion paper: the theory establishes the variational mechanism \cite{MendicoTheory2026}, whereas the present work demonstrates that the mechanism yields biologically interpretable, anatomically constrained and dynamically non-trivial reaction maps in a multimodal setting.

\section*{Results}

\subsection*{Task-evoked fMRI defines spatially resolved stimulation and reaction signatures}

\begin{figure*}[!t]
    \centering
    \includegraphics[width=\textwidth]{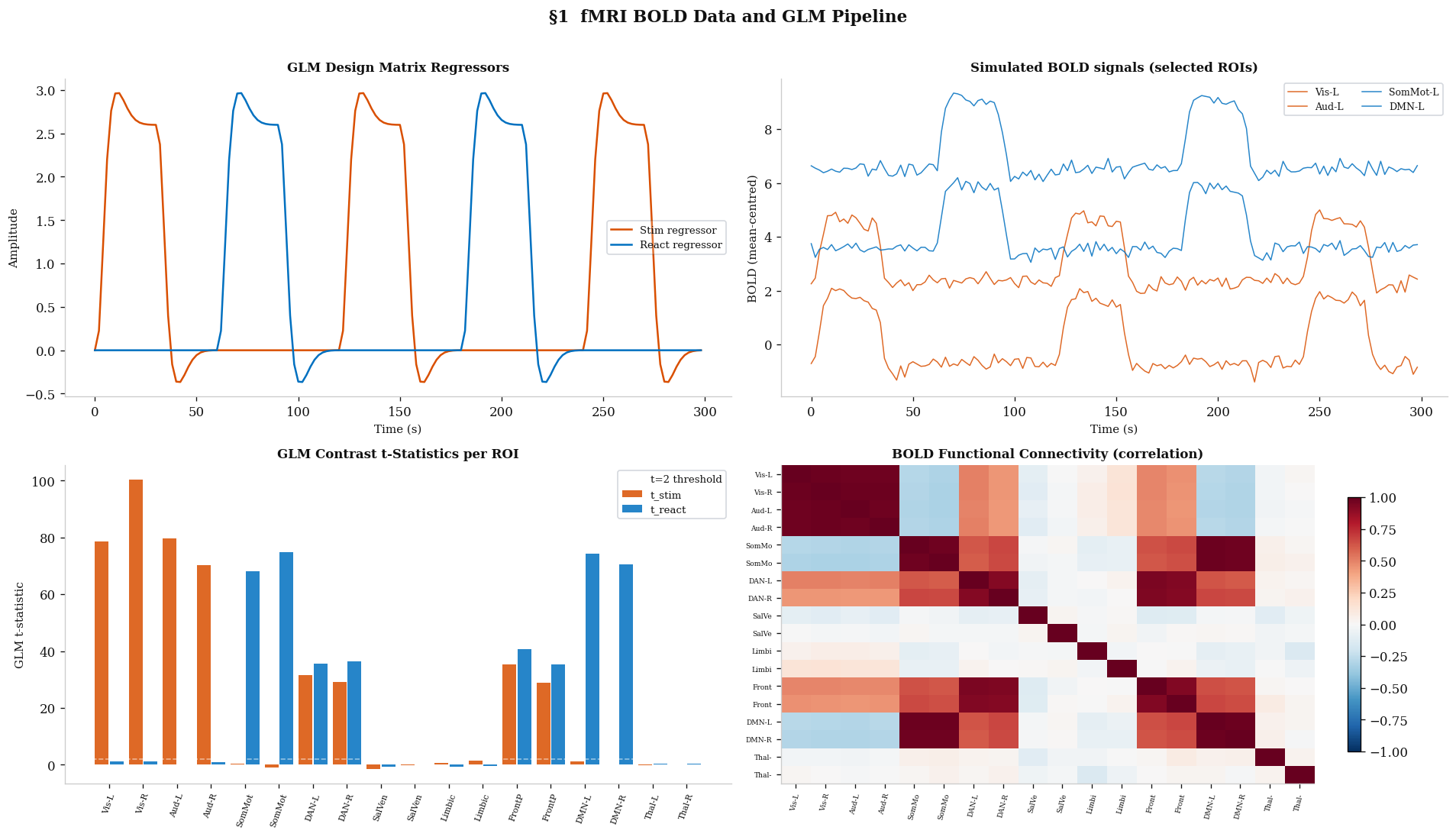}
    \caption{\textbf{Task-evoked fMRI BOLD pipeline.}
    Simulated blood-oxygen-level-dependent data analysed with a general linear model. Left, stimulus and reaction regressors after convolution with the haemodynamic response. Top right, representative blood-oxygen-level-dependent signals in selected regions of interest. Bottom left, region-wise contrast statistics showing strong task selectivity. Bottom right, blood-oxygen-level-dependent functional connectivity. This figure defines the fMRI component of the source and target estimation step.}
    \label{fig:fmri}
\end{figure*}

We first simulated a task-based blood-oxygen-level-dependent experiment on a common cortical support with $18$ regions of interest. As shown in Fig.~\ref{fig:fmri}, the block-design general linear model cleanly separates a stimulation regressor from a reaction regressor and generates region-wise contrast statistics with strong spatial selectivity. Early task epochs concentrate on visually and auditorily driven regions, whereas later task epochs emphasise sensorimotor and default-mode-related regions. The resulting blood-oxygen-level-dependent correlation structure is non-trivial and already exhibits mesoscale organisation rather than independent regional responses.

These outputs are important for two reasons. First, they show that the source and target measures can be grounded in standard task-based fMRI statistics rather than hand-crafted weights \cite{Friston2009,Esteban2020}. Second, they provide a spatially localised but temporally coarse characterisation of the stimulus-to-reaction transformation, which is precisely the role that fMRI should play in the multimodal pipeline.

\subsection*{Source-reconstructed EEG/MEG resolves the temporal separation between stimulus and reaction}

\begin{figure*}[!t]
    \centering
    \includegraphics[width=\textwidth]{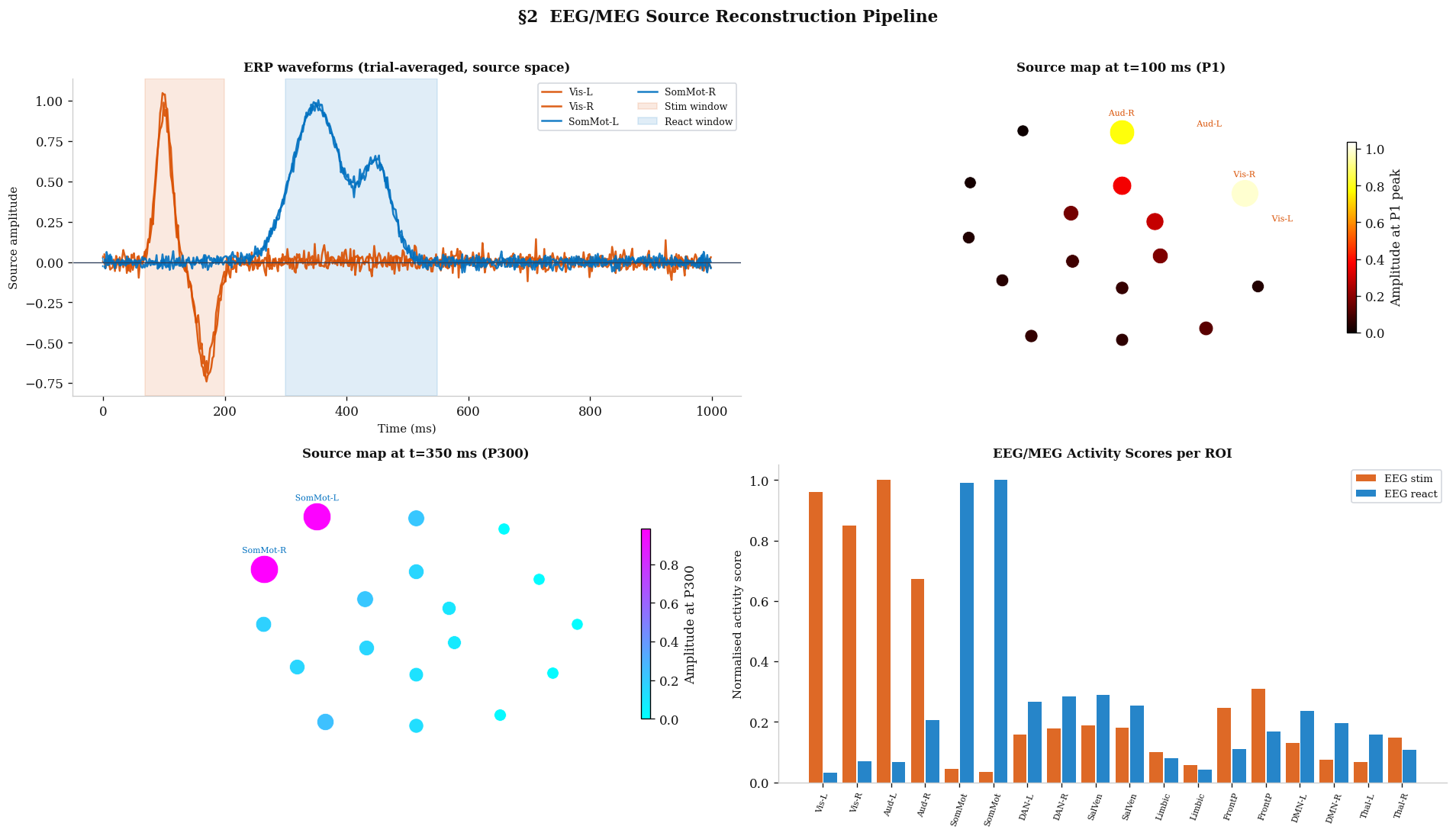}
    \caption{\textbf{Source-reconstructed EEG/MEG pipeline.}
    Top left, trial-averaged source-space ERP waveforms showing early stimulus-locked and later reaction-locked components. Top right, source map at $\sim 100$ ms. Bottom left, source map at $\sim 350$ ms. Bottom right, region-wise electrophysiological activity scores for the stimulation and reaction windows. This figure provides the temporally resolved component of the source and target estimation step.}
    \label{fig:eeg}
\end{figure*}

We next constructed source-reconstructed electrophysiological data using a lead-field forward model and a minimum-norm inverse procedure. Figure~\ref{fig:eeg} shows that this modality complements blood-oxygen-level-dependent data exactly in the way needed for the transport problem: early components concentrate in sensory entry regions, whereas later components dominate in reaction-related regions. The spatial source maps at approximately $100$ ms and $350$ ms reveal a clean temporal shift from stimulus-locked to reaction-locked activity.

This temporal separation is crucial. The source measure $\mu_{\mathrm{stim}}^{+}$ should reflect entry of the external stimulation into the effective propagation architecture, whereas $\mu_{\mathrm{react}}^{-}$ should reflect the later reaction-producing configuration. Blood-oxygen-level-dependent data alone localise these patterns but blur their timing; source-reconstructed EEG/MEG resolves them sharply. In this sense, Fig.~\ref{fig:eeg} provides the temporal half of the inverse problem. The source-reconstruction viewpoint is standard in modern EEG/MEG analysis and is particularly useful when one aims to distinguish early stimulus-driven components from later reaction-related components \cite{Knosche2022}.

\newpage

\subsection*{Multimodal fusion yields balanced source and target probability measures}

We then combined the spatial specificity of blood-oxygen-level-dependent responses and the temporal precision of source-reconstructed electrophysiology through weighted geometric fusion of modality-specific activity scores. If $a_i^{\mathrm{fMRI}}$ and $a_i^{\mathrm{EEG}}$ denote normalised regional scores, the fused stimulation profile is
\[
a_i^{\mathrm{stim}}
=
\bigl(a_i^{\mathrm{fMRI,stim}}\bigr)^{w_f}
\bigl(a_i^{\mathrm{EEG,stim}}\bigr)^{1-w_f},
\]
with $w_f \in [0,1]$, and similarly for the reaction profile. After normalisation we obtain the balanced measures
\[
\mu_{\mathrm{stim}}^{+},
\qquad
\mu_{\mathrm{react}}^{-},
\qquad
b(v)=\mu_{\mathrm{stim}}^{+}-\mu_{\mathrm{react}}^{-}.
\]

\begin{figure*}[!t]
    \centering
    \includegraphics[width=\textwidth]{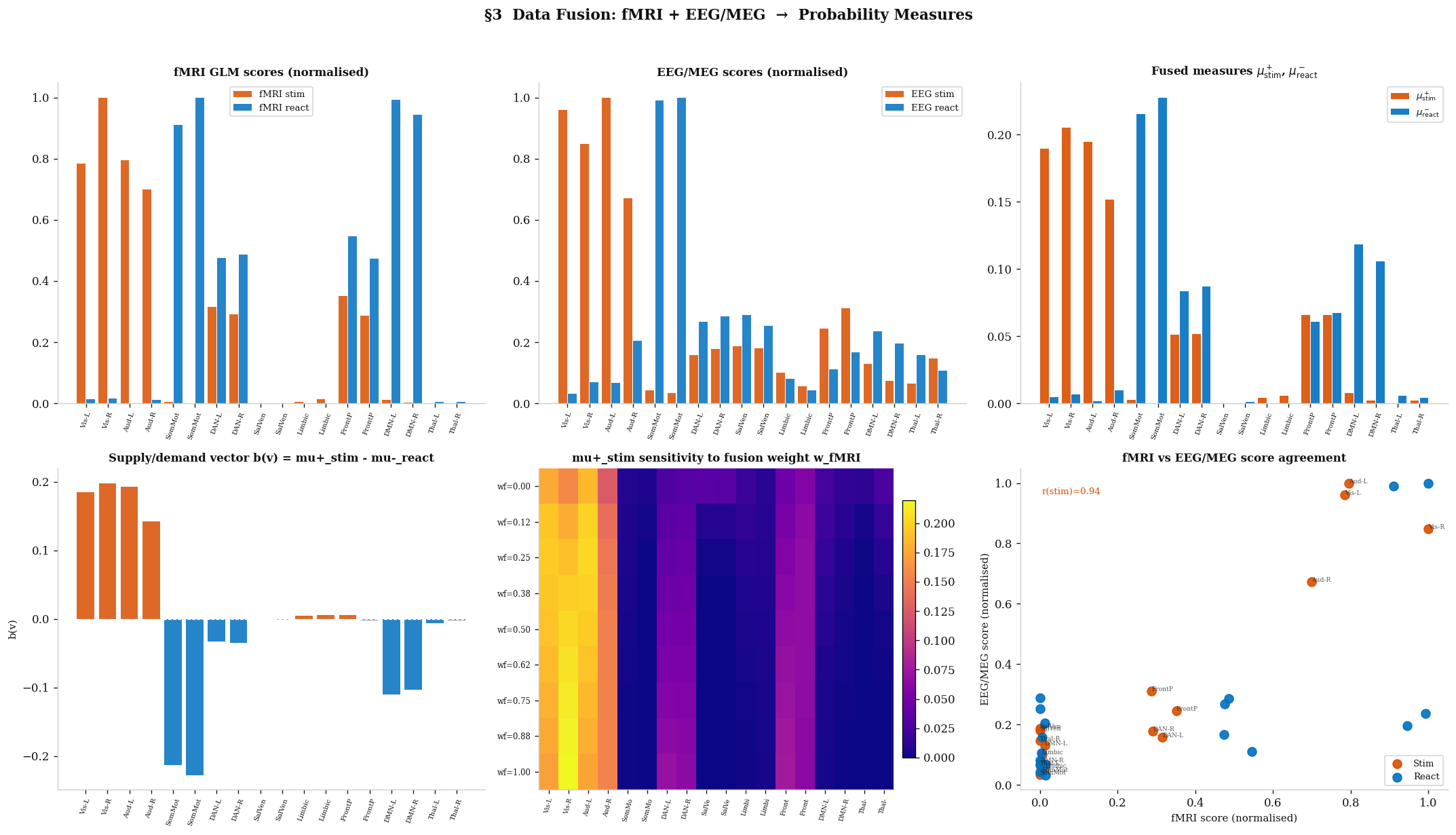}
    \caption{\textbf{Multimodal fusion and probability measures.}
    Top row, normalised fMRI scores, normalised EEG/MEG scores, and fused source and target measures. Bottom left, supply--demand vector $b(v)=\mu_{\mathrm{stim}}^{+}-\mu_{\mathrm{react}}^{-}$. Bottom middle, sensitivity of $\mu_{\mathrm{stim}}^{+}$ to the fusion weight. Bottom right, agreement between blood-oxygen-level-dependent and electrophysiological regional scores. This figure defines the balanced source and target measures used in the branched transport problem.}
    \label{fig:fusion}
\end{figure*}

Figure~\ref{fig:fusion} shows the resulting modality-specific patterns and fused measures. The fused $\mu_{\mathrm{stim}}^{+}$ remains concentrated on sensory entry regions, while the fused $\mu_{\mathrm{react}}^{-}$ shifts toward motor and higher-order reaction regions. At the same time, intermediate relay regions acquire non-zero mass through cross-modal agreement rather than arbitrary assignment.

The lower row of Fig.~\ref{fig:fusion} shows two additional points. First, the supply--demand vector $b(v)$ is sharply structured, producing a transport problem that is neither trivial nor diffuse. Second, the sensitivity map with respect to the fusion parameter demonstrates that the dominant support is stable over a broad range of modality weightings, although quantitative mass redistribution persists in secondary regions. The transport problem is therefore not being driven by a fragile or modality-specific artefact. Additional analyses across the full \(\alpha\)-grid confirm that the transition from weakly aggregated to strongly aggregated routing is continuous over the stable parameter regime (Supplementary Figs.~S1 and S2).

\subsection*{Diffusion-informed anisotropy induces a directional transport geometry}

To encode anatomical directionality, we assigned to each region of interest a synthetic diffusion tensor and used it to build a local anisotropic cost
\[
c(x,\tau)=\sqrt{\tau^\top (D(x)+\varepsilon I)^{-1}\tau}.
\]
For each candidate arc \(e\), this induces a graph-level cost
\[
\beta_e
=
\int_0^{\ell_e}
c\!\left(\gamma_e(s),\frac{\dot{\gamma}_e(s)}{\|\dot{\gamma}_e(s)\|}\right)\,ds,
\]
approximated numerically by midpoint quadrature along the edge.

\begin{figure*}[!t]
    \centering
    \includegraphics[width=\textwidth]{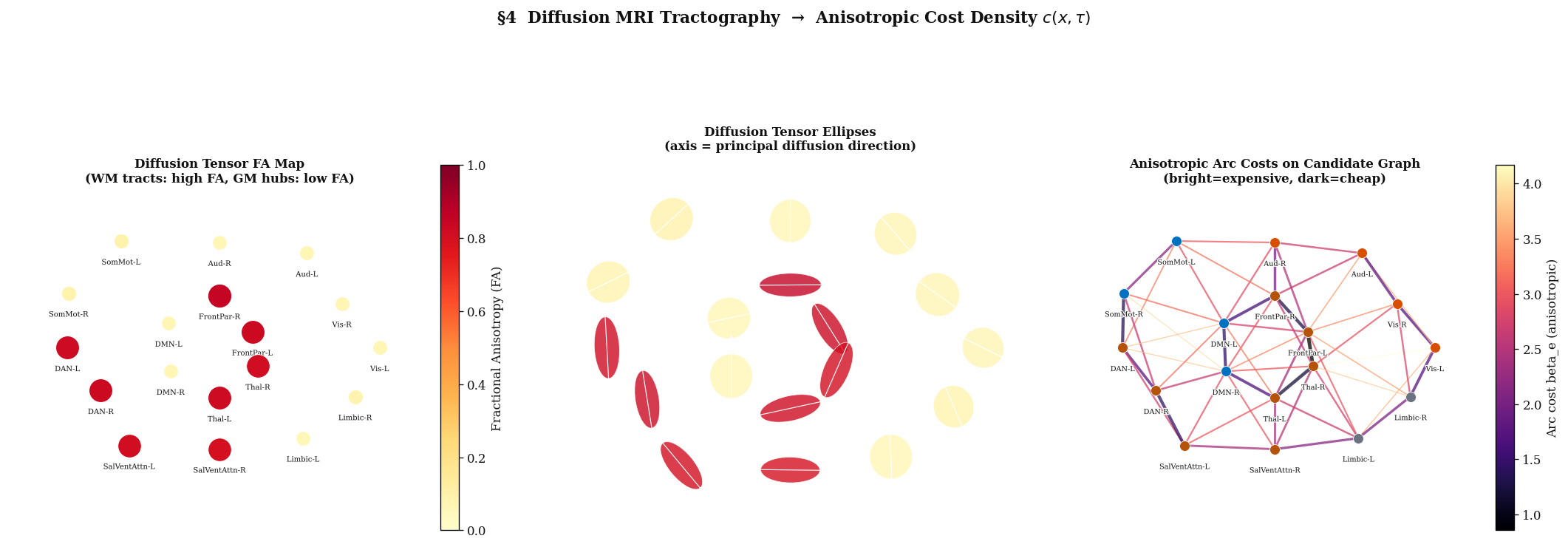}
    \caption{\textbf{Diffusion-informed anisotropic transport geometry.}
    Left, synthetic fractional-anisotropy map across regions of interest. Middle, diffusion tensor ellipses and principal diffusion axes. Right, anisotropic arc costs on the candidate graph, showing strong directional heterogeneity. This figure defines the tractography-informed transport cost used in the anisotropic optimisation.}
    \label{fig:dti}
\end{figure*}

Figure~\ref{fig:dti} shows the resulting tractography-informed prior: white-matter-like relay regions have high fractional anisotropy, principal diffusion axes introduce directional preferences, and graph-level arc costs become heterogeneous even among geometrically similar connections. This construction is directly inspired by the use of diffusion-informed structure to constrain large-scale communication models in connectomics \cite{Hagmann2007,Yeh2021,Zhang2022}.

This step is decisive because it transforms the transport problem from a distance-based optimisation into a geometry-aware inverse problem. In the isotropic baseline, short edges are generically favoured. In the anisotropic model, some short edges become expensive if they cut across implausible directions, while longer edges can become favourable if they align with dominant tensor axes. Thus, the anatomical prior does not simply modulate the cost magnitude; it changes the directional logic of admissible propagation.

\subsection*{Anisotropic branched transport infers a distinct reaction-map backbone}

We then solved the discrete branched transport problem
\[
\min_{w\ge 0}\;
\sum_{e\in\mathcal E}\beta_e\,w_e^\alpha
\qquad \text{subject to} \qquad
Aw=b,
\]
with \(0<\alpha<1\), where the concavity of \(w\mapsto w^\alpha\) promotes aggregation and branching \cite{Xia2003,MendicoTheory2026}.

\begin{figure*}[!t]
    \centering
    \includegraphics[width=\textwidth]{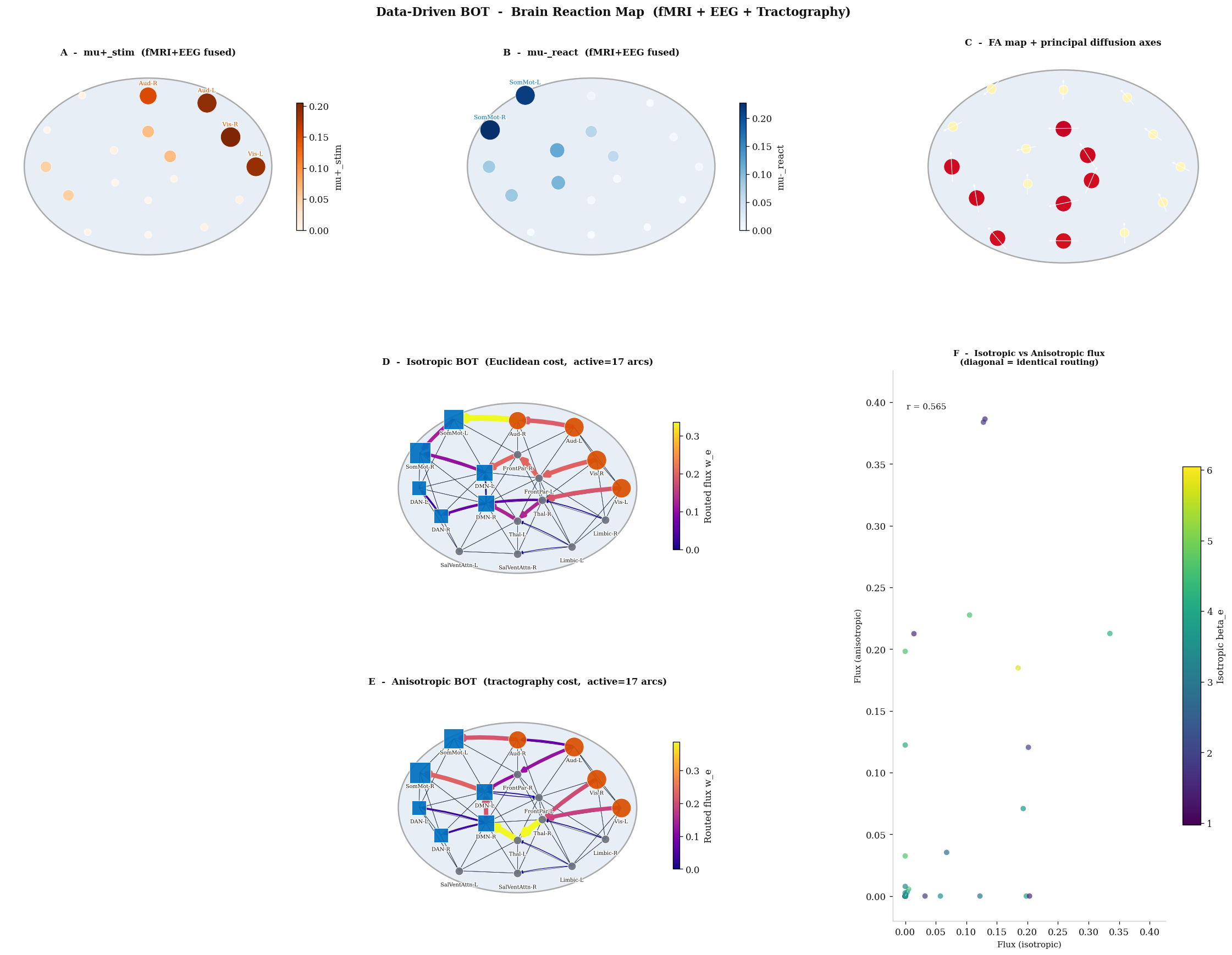}
    \caption{\textbf{Data-driven branched optimal transport infers the brain reaction map.}
    (a) Fused $\mu_{\mathrm{stim}}^{+}$. (b) Fused $\mu_{\mathrm{react}}^{-}$. (c) Fractional-anisotropy map and principal diffusion axes. (d) Isotropic branched transport solution. (e) Anisotropic branched transport solution. (f) Edgewise flux comparison between isotropic and anisotropic models. The anisotropic model yields a routing backbone that is qualitatively distinct from the isotropic baseline and more strongly aligned with anatomical directional priors.}
    \label{fig:botmap}
\end{figure*}

Figure~\ref{fig:botmap} is the central result of the paper. The isotropic solution yields a branched architecture but remains relatively close to direct geometrical routing. The anisotropic solution, by contrast, reorganises the routing backbone around relay regions aligned with the tractography-derived geometry. Several branches that are only weakly expressed in the isotropic case become dominant under anisotropic costs, while some direct alternatives disappear entirely.

The flux comparison in Fig.~\ref{fig:botmap}f shows that the two solutions differ not only by small perturbations in weight but by substantial redistributions of mass across edges. In other words, anisotropy qualitatively reshapes the inferred reaction map. This is exactly the type of effect that fixed-substrate control models cannot reveal: once the substrate is prescribed, anatomical anisotropy can modulate dynamics on that substrate, but it cannot alter which routing architecture is selected as the canonical backbone in the first place. Comparison with simpler non-branched baselines further shows that the concave transport cost is specifically responsible for the emergence of shared relay corridors (Supplementary Fig.~S3).

\subsection*{Hybrid stochastic dynamics reveals geometric--dynamical trade-offs and rank reversals}

Finally, we used the anisotropic optimal graph as the substrate of the hybrid stochastic extension. On the inferred graph we defined a linear stochastic dynamics of the form
\[
dX_t
=
A_{G,w}X_t\,dt
+
B_{\mathrm{stim}}a_t\,dt
+
u_t\,dt
+
C_{G,w}\,dW_t,
\]
with
\[
A_{G,w}=-\kappa I-\beta_{\mathrm{dyn}}L_{G,w},
\qquad
C_{G,w}=\sigma_0 I+\sigma_1 D_{G,w}^{1/2}.
\]
The associated hybrid functional is
\[
F_\lambda(G,w)=E_\alpha(G,w)+\lambda J_{\mathrm{dyn}}(G,w).
\]

\begin{figure*}[!t]
    \centering
    \includegraphics[width=\textwidth]{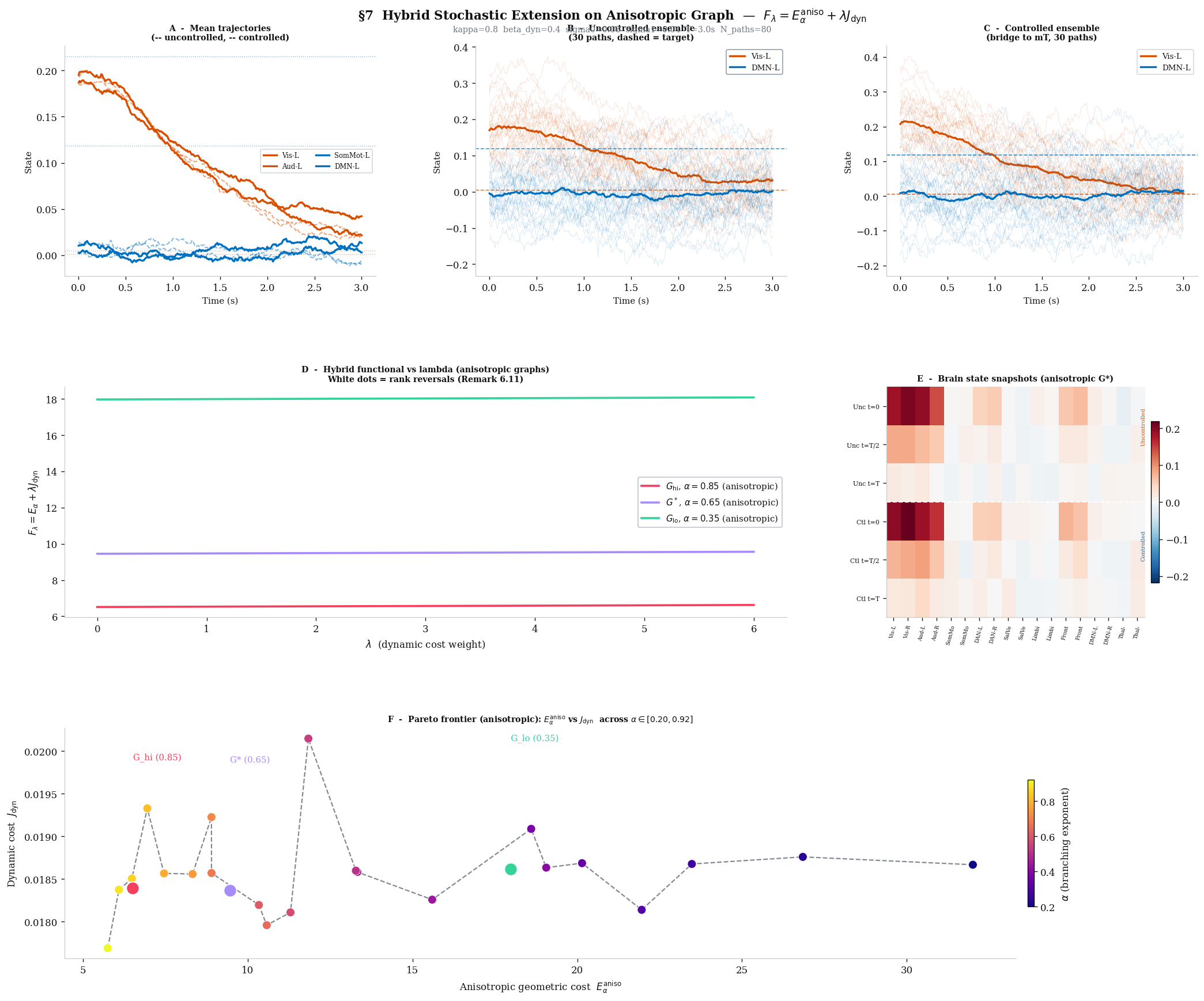}
    \caption{\textbf{Hybrid stochastic extension on the anisotropic graph.}
    (a) Mean trajectories under uncontrolled and controlled dynamics. (b) Uncontrolled path ensemble. (c) Controlled path ensemble. (d) Hybrid functional $F_\lambda$ for competing graphs, showing rank reversals as $\lambda$ varies. (e) Brain-state snapshots. (f) Pareto frontier across branching exponents $\alpha$. This figure shows that geometric efficiency and dynamical controllability define distinct and competing criteria for reaction-map selection.}
    \label{fig:hybrid}
\end{figure*}

Figure~\ref{fig:hybrid} shows three complementary results. First, controlled stochastic trajectories reach terminal states substantially closer to the prescribed reaction profile than uncontrolled trajectories, confirming that the inferred graph is dynamically usable as a state-transition substrate. Second, the hybrid functional reveals that graphs with lower geometric cost are not necessarily dynamically preferred. Third, the Pareto frontier across branching exponents $\alpha$ exhibits a non-trivial geometry, with local trade-offs between ramification and dynamic cost.

The geometric meaning of Fig.~\ref{fig:hybrid}f is intuitive: each point corresponds to one candidate routing map, moving left reduces geometric transport cost, and moving down reduces dynamic control cost. The Pareto frontier is the set of non-dominated maps for which one criterion cannot be improved without worsening the other. Most importantly, Fig.~\ref{fig:hybrid}d shows rank reversals across $\lambda$. As the weight of the dynamical term increases, the ordering of candidate graphs changes. A graph that is geometrically optimal can lose to a dynamically cheaper competitor, and vice versa. This is a strong argument against using purely geometric or purely dynamical criteria in isolation. The natural object is the coupled geometric--dynamical landscape. The graph-induced control viewpoint is consistent with recent work on stochastic steering and Schrödinger bridge formulations for brain state transitions \cite{Kawakita2022,Kamiya2023}.

\section*{Discussion}

The main outcome of this study is that multimodal neuroimaging, tractography-informed anisotropy and branched transport can be integrated into a single inverse framework that infers propagation architecture rather than presupposing it. The inferred object is not a trajectory, not a control signal and not a predefined connectome, but a \emph{brain reaction map}: an effective routing backbone that connects stimulation to distributed response under anatomical and dynamical constraints.

Three conclusions emerge. First, multimodal estimation of source and target measures is both feasible and structurally informative. The combination of blood-oxygen-level-dependent spatial localisation and electrophysiological temporal resolution produces source and target measures that are sufficiently concentrated to drive a non-trivial transport problem while remaining robust to the fusion weight. Second, anisotropy matters at the architectural level. The anisotropic branched solution is not simply a weighted perturbation of the isotropic one; it activates a different routing logic, with stronger relay structure and clearer anatomical alignment. Third, dynamic plausibility does not reduce to geometric efficiency. The hybrid extension shows that the graphs preferred by transport cost alone need not minimise the dynamic control cost, and that rank reversals across $\lambda$ are intrinsic rather than exceptional.

The biological interpretability of the inferred architecture is an important point. In the present synthetic setting, the relay regions selected by the anisotropic BOT solution are not anonymous intermediates: the strongest bottlenecks and branching interfaces emerge preferentially in dorsal attention, salience/ventral attention, frontoparietal and thalamic nodes, that is, in systems commonly interpreted as integrative bridges between sensory input and distributed action or cognitive output \cite{Bazinet2023,Fotiadis2024,Gilson2020,Patow2024}. In this sense, the model does not simply recover short routes between sources and sinks; it selects a mesoscale backbone in which anatomically and functionally plausible association systems act as shared highways for signal aggregation before redistribution. Additional discussion of the mesoscale interpretation of relay systems is provided in Supplementary Note 5.

The interpretation of the Pareto frontier is also important for readers outside optimal transport or variational analysis. Each point in the \((E_\alpha,J_{\mathrm{dyn}})\) plane represents one candidate routing architecture, with one axis measuring geometric economy and the other dynamic steering cost. The frontier therefore marks the set of graphs that are optimal trade-offs between these two criteria. Rank reversals occur when the preferred architecture changes as the relative weight of dynamics is increased: a graph may be the most economical one geometrically, yet cease to be optimal once controllability is taken seriously. This is not a technical curiosity but a conceptual message: geometry and dynamics should be treated as coupled principles of large-scale propagation rather than as interchangeable surrogates.

Several limitations remain. The present article uses synthetic multimodal data and a low-dimensional cortical support with \(18\) regions of interest, the candidate graph is finite, and the stochastic dynamics is linear. The tractography prior is represented by a coarse tensor field rather than subject-specific whole-brain tract reconstruction. These simplifications are intentional: they isolate the effect of multimodal source--target inference, anisotropic geometry and hybrid selection in a transparent setting. The next step is therefore clear: the framework should be tested on real task-based fMRI, source-localised EEG/MEG and tractography-derived costs.

A related concern is scalability. The current implementation solves a constrained nonlinear optimisation problem with SLSQP on a sparse candidate graph and then estimates dynamic costs by Monte Carlo simulation. For the present proof-of-concept, this strategy is stable and sufficient. Scaling to human connectomes with 200 or more regions is conceptually straightforward but computationally more demanding, because the number of candidate arcs, optimisation variables and repeated dynamic evaluations all increase substantially. In practice, this points to several natural extensions: anatomically informed graph sparsification, continuation or warm-start strategies across \(\alpha\), parallel evaluation of dynamic costs, and alternative large-scale optimisation schemes based on decomposition, proximal updates or differentiable surrogates. The current results should therefore be read not as a final engineering solution for whole-brain resolution, but as a variational proof-of-principle that motivates these larger-scale developments. Practical routes toward scaling the framework beyond the present 18-ROI proof-of-concept are discussed in Supplementary Note 6 and Supplementary Fig.~S4.

More broadly, the framework suggests that a substantial part of what is often interpreted as state-transition difficulty may in fact be a question of architecture inference. If the propagation substrate itself is allowed to vary, then the correct inverse problem is no longer only to steer dynamics on a graph, but to determine which graph best explains how stimulation is transformed into reaction. In this sense, the main computational contribution of the present article should be read together with the theoretical foundation developed in the companion paper \cite{MendicoTheory2026}.

\section*{Methods}

\subsection*{Overview of the computational pipeline}

The computational workflow consists of five successive stages: simulation of task-related blood-oxygen-level-dependent data, simulation of source-reconstructed electrophysiology, multimodal fusion of stimulation and reaction measures, construction of anisotropic transport costs from diffusion-informed tensors, and inference of a reaction map through branched transport optimisation followed by graph-induced stochastic analysis. All simulations were implemented in Python in a companion notebook.

\subsection*{Synthetic cortical support and region layout}

All computations were performed on a common support made of \(N_{\mathrm{roi}}=18\) regions of interest embedded in a two-dimensional cortical slice. Twelve regions were placed on an outer elliptical shell and six on an inner shell, producing a coarse geometry with peripheral sensory-association regions and inner relay regions. Regions were labelled to mimic visual, auditory, sensorimotor, dorsal attention, salience/ventral attention, limbic, frontoparietal, default mode and thalamic systems.

\subsection*{Task-based fMRI simulation and general linear model}

The fMRI acquisition was simulated over \(T_{\mathrm{fmri}}=300\) s with repetition time \(TR=2.0\) s, yielding \(N_{\mathrm{vol}}=150\) volumes. Stimulus and reaction epochs followed a repeating block design consisting of 30 s stimulus, 30 s rest, 30 s reaction and 30 s rest. Neural task regressors were convolved with a canonical double-gamma haemodynamic response function
\[
h(t)=
\frac{t^{a_1-1}e^{-t/b_1}}{\Gamma(a_1)b_1^{a_1}}
-
c\,\frac{t^{a_2-1}e^{-t/b_2}}{\Gamma(a_2)b_2^{a_2}},
\qquad t\ge 0,
\]
with \(a_1=6\), \(a_2=16\), \(b_1=b_2=1\), \(c=1/6\), sampled on \([0,32]\) s and normalised to unit peak magnitude.

The design matrix was
\[
X_{\mathrm{GLM}}=
\bigl[x_{\mathrm{stim}},x_{\mathrm{react}},\mathbf{1}\bigr],
\]
where \(x_{\mathrm{stim}}\) and \(x_{\mathrm{react}}\) denote the haemodynamically convolved task regressors and \(\mathbf{1}\) is the intercept. Region-wise blood-oxygen-level-dependent data were generated as
\[
Y=X_{\mathrm{GLM}}\beta_{\mathrm{true}}+\varepsilon,
\]
with Gaussian noise of standard deviation \(0.15\). Positive stimulus- and reaction-related contrast statistics were extracted as the fMRI activity scores.

\subsection*{Source-reconstructed EEG/MEG simulation}

Source-reconstructed electrophysiological data were simulated over a \(T_{\mathrm{eeg}}=1.0\) s epoch with sampling rate \(f_s=512\) Hz, giving \(N_{\mathrm{samp}}=512\) samples per epoch. Sensor-level data were generated on \(N_{\mathrm{sensor}}=64\) sensors placed on a circle around the cortical support. The lead field \(L\in\mathbb{R}^{N_{\mathrm{sensor}}\times N_{\mathrm{roi}}}\) was defined by inverse-square distance falloff between sensors and source regions and then normalised by its maximum entry.

Source reconstruction used a minimum-norm inverse
\[
L^\dagger = L^\top(LL^\top+\lambda_{\mathrm{reg}}I)^{-1},
\]
with Tikhonov regularisation parameter \(\lambda_{\mathrm{reg}}=0.05\). Early stimulus-locked source components were placed in visual and auditory regions, whereas later reaction-locked components were placed in sensorimotor and default-mode-related regions. Forty noisy trials were generated and averaged to improve signal-to-noise ratio. Stimulation and reaction source scores were obtained from mean absolute source amplitudes over the windows \([70,200]\) ms and \([300,550]\) ms, respectively.

\subsection*{Multimodal fusion of source and target measures}

Let \(a_i^{\mathrm{fMRI,stim}}\), \(a_i^{\mathrm{fMRI,react}}\), \(a_i^{\mathrm{EEG,stim}}\), and \(a_i^{\mathrm{EEG,react}}\) denote the modality-specific regional scores after normalisation to unit maximum. The fused stimulation and reaction profiles were defined by weighted geometric averaging,
\[
a_i^{\mathrm{stim}}
=
\bigl(a_i^{\mathrm{fMRI,stim}}+\varepsilon_0\bigr)^{w_f}
\bigl(a_i^{\mathrm{EEG,stim}}+\varepsilon_0\bigr)^{1-w_f},
\]
\[
a_i^{\mathrm{react}}
=
\bigl(a_i^{\mathrm{fMRI,react}}+\varepsilon_0\bigr)^{w_f}
\bigl(a_i^{\mathrm{EEG,react}}+\varepsilon_0\bigr)^{1-w_f},
\]
with \(w_f=0.55\) and \(\varepsilon_0=10^{-6}\). The final source and target measures were
\[
\mu_{\mathrm{stim}}^{+}(i)=\frac{a_i^{\mathrm{stim}}}{\sum_j a_j^{\mathrm{stim}}},
\qquad
\mu_{\mathrm{react}}^{-}(i)=\frac{a_i^{\mathrm{react}}}{\sum_j a_j^{\mathrm{react}}},
\]
and the supply--demand vector was \(b=\mu_{\mathrm{stim}}^{+}-\mu_{\mathrm{react}}^{-}\).

Sensitivity of the fused stimulation measure to the fusion weight was assessed over a grid \(w_f\in\{0,0.125,\dots,1\}\).

\subsection*{Diffusion-informed anisotropic cost construction}

Synthetic diffusion tensors \(D_i\in\mathbb{R}^{2\times2}\) were assigned to each region of interest. Relay regions corresponding to dorsal attention, salience/ventral attention, frontoparietal and thalamic systems were given higher anisotropy, whereas other regions were assigned lower anisotropy. Fractional anisotropy was computed from the eigenvalues of each tensor.

A candidate graph was built by a \(k\)-nearest-neighbour rule with \(k=5\), producing an undirected edge set that was duplicated into directed arcs. For an arc \(e=(i,j)\) with tangent direction \(\tau_{ij}\) and Euclidean length \(\ell_{ij}\), the anisotropic edge cost was defined by midpoint evaluation:
\[
\beta_{ij}^{\mathrm{aniso}}
=
\ell_{ij}
\sqrt{
\tau_{ij}^{\top}
\bigl(D_{ij}^{\mathrm{mid}}+\varepsilon I\bigr)^{-1}
\tau_{ij}
},
\]
where \(D_{ij}^{\mathrm{mid}}=\frac12(D_i+D_j)\) and \(\varepsilon=0.05\). For comparison, an isotropic baseline cost was also constructed from Euclidean distance and a scalar white-matter score.

\subsection*{Discrete branched transport optimisation}

Let \(A\) denote the incidence matrix of the directed candidate graph. For a fixed branching exponent \(\alpha\in(0,1)\), the optimisation problem was
\[
\min_{w\ge 0}\;
\sum_{e\in\mathcal E}\beta_e\,w_e^\alpha
\qquad \text{subject to}\qquad
Aw=b.
\]
The principal analyses used \(\alpha=0.65\). Optimisation was performed with SLSQP and 12 random restarts. Each restart was initialised from a combination of a clipped least-squares solution of the balance constraints and a small positive random perturbation. After optimisation, fluxes below \(10^{-4}\) times the maximal flux were set to zero for support visualisation only. The same candidate graph and optimisation procedure were used for both isotropic and anisotropic costs.

\subsection*{Graph-induced stochastic dynamics and dynamic cost}

The anisotropic optimal graph was converted into a directed weighted adjacency matrix \(W\), then symmetrised to \(S=\frac12(W+W^\top)\). The weighted degree vector \(d\) and graph Laplacian \(L\) were defined by
\[
d_i=\sum_j S_{ij},
\qquad
L=\mathrm{diag}(d)-S.
\]
The graph-induced stochastic dynamics was
\[
dX_t
=
A_{G,w}X_t\,dt
+
B_{\mathrm{stim}}a_t\,dt
+
u_t\,dt
+
C_{G,w}\,dW_t,
\]
with
\[
A_{G,w}=-\kappa I-\beta_{\mathrm{dyn}}L,
\qquad
C_{G,w}=\sigma_0 I+\sigma_1 \mathrm{diag}(\sqrt{d}).
\]
The parameters were \(\kappa=0.8\), \(\beta_{\mathrm{dyn}}=0.4\), \(\sigma_0=0.08\), and \(\sigma_1=0.04\). The simulation horizon was \(T_{\mathrm{sim}}=3.0\) s with timestep \(dt=0.005\), giving \(N_{\mathrm{steps}}=600\) steps. Path ensembles contained \(N_{\mathrm{paths}}=80\) trajectories.

The uncontrolled process corresponded to \(u_t\equiv 0\). The controlled process used a bridge-type feedback toward the target profile,
\[
u_t
=
C_{G,w}C_{G,w}^{\top}
\frac{m_T-X_t}{T_{\mathrm{sim}}-t+\varepsilon_{\mathrm{bridge}}},
\]
with \(\varepsilon_{\mathrm{bridge}}=0.05\). Dynamic cost was estimated by Monte Carlo approximation of the quadratic control functional associated with the controlled process.

\subsection*{Pareto and \texorpdfstring{$\lambda$}{lambda}-analysis}

The anisotropic transport problem was repeated over a grid of branching exponents
\[
\alpha \in \{0.20,0.20+\Delta,\dots,0.92\},
\]
with 22 evenly spaced values. For each inferred graph, the dynamic cost \(J_{\mathrm{dyn}}\) was estimated numerically. The hybrid functional
\[
F_\lambda(G,w)=E_\alpha(G,w)+\lambda J_{\mathrm{dyn}}(G,w)
\]
was evaluated over \(\lambda\in[0,6]\) on a grid of 300 values. Pareto points were obtained from the set of \((E_\alpha,J_{\mathrm{dyn}})\) pairs and rank reversals were identified by sign changes in pairwise differences of the corresponding \(F_\lambda\) curves.

\section*{Data availability}

This study uses synthetic multimodal datasets generated within the companion computational notebook. All scripts required to reproduce the simulations, figures and derived measures will be made publicly available upon publication.

\section*{Code availability}

The full pipeline, including fMRI simulation, source-reconstructed EEG/MEG, multimodal fusion, anisotropic branched transport optimisation and stochastic graph dynamics, will be deposited in a public repository upon publication.





\section*{Competing interests}

The author declares no competing interests.

\newpage

\appendix

\section*{Supplementary overview}

This Supplementary Information is designed to reinforce the central claims of the main manuscript along four dimensions that are especially relevant for evaluation by a broad multidisciplinary readership. First, it examines how the inferred anisotropic routing architecture varies with the branching exponent \(\alpha\), thereby clarifying that the reported reaction maps do not arise from a single isolated parameter choice. Second, it compares the anisotropic branched solution with simpler non-branched baselines, making explicit which features of the inferred architecture are genuinely due to the ramified objective. Third, it quantifies the relay-region structure of the inferred map, thereby strengthening the biological interpretability of the selected transport backbone. Fourth, it provides a proof-of-principle scalability analysis, addressing the natural question of whether the present framework can plausibly extend beyond the \(18\)-ROI synthetic demonstration.

Taken together, these supplementary analyses support four stronger versions of the claims made in the main text. The first is that the anisotropic reaction map is \emph{robust across a stable branching regime}, rather than being a fragile solution at one preferred \(\alpha\). The second is that the emergence of \emph{shared relay corridors} is a distinctive signature of branched transport and is not reproduced by simpler linear-cost or shortest-path surrogates. The third is that the strongest relay regions are \emph{biologically interpretable mesoscale bottlenecks}, not anonymous graph intermediates. The fourth is that the current implementation should be read as a \emph{variational proof-of-principle with plausible computational paths to scale-up}, rather than as an already final large-connectome engineering solution.

\section*{Supplementary Note 1. Extended dependence on the branching exponent \texorpdfstring{$\alpha$}{alpha}}

A defining ingredient of the model is the concavity of the edge cost \(w_e^\alpha\), with \(0<\alpha<1\). Smaller values of \(\alpha\) favour stronger route sharing and earlier aggregation of flux; values closer to \(1\) weaken the ramification incentive and produce more distributed routing. In the main manuscript, the principal reaction maps are shown for a representative value \(\alpha=0.65\). Here we make explicit how the global observables and inferred supports behave across the broader explored regime.

\subsection*{Global observables across the \texorpdfstring{$\alpha$}{alpha}-grid}

Supplementary Fig.~S1 reports the anisotropic geometric cost \(E_\alpha\), the dynamic cost \(J_{\mathrm{dyn}}\), and the support size \(|e^\ast|\) over the full \(\alpha\)-grid. The dominant trend away from the transition region is clear: the geometric cost decreases as \(\alpha\) increases, while the support size remains in a relatively stable range. The dynamic cost varies much more weakly over the stable regime than the geometric cost, which is one of the reasons why the Pareto structure in the main manuscript is non-trivial.

A particularly important feature is the sharp anomaly around \(\alpha \approx 0.4\), where both the geometric and dynamic observables undergo abrupt changes and the support collapses. We interpret this interval as a \emph{near-degenerate transition regime} of the optimisation landscape rather than as the biologically relevant operating point of the framework. This regime is scientifically informative: it indicates that the model separates naturally into weakly branched and strongly branched phases, with a narrow transition region between them. However, the main conclusions of the article do not rely on that transition itself, but on the stable family of solutions around and above the main-text reference value \(\alpha=0.65\).

\begin{figure*}[!t]
    \centering
    \includegraphics[width=\textwidth]{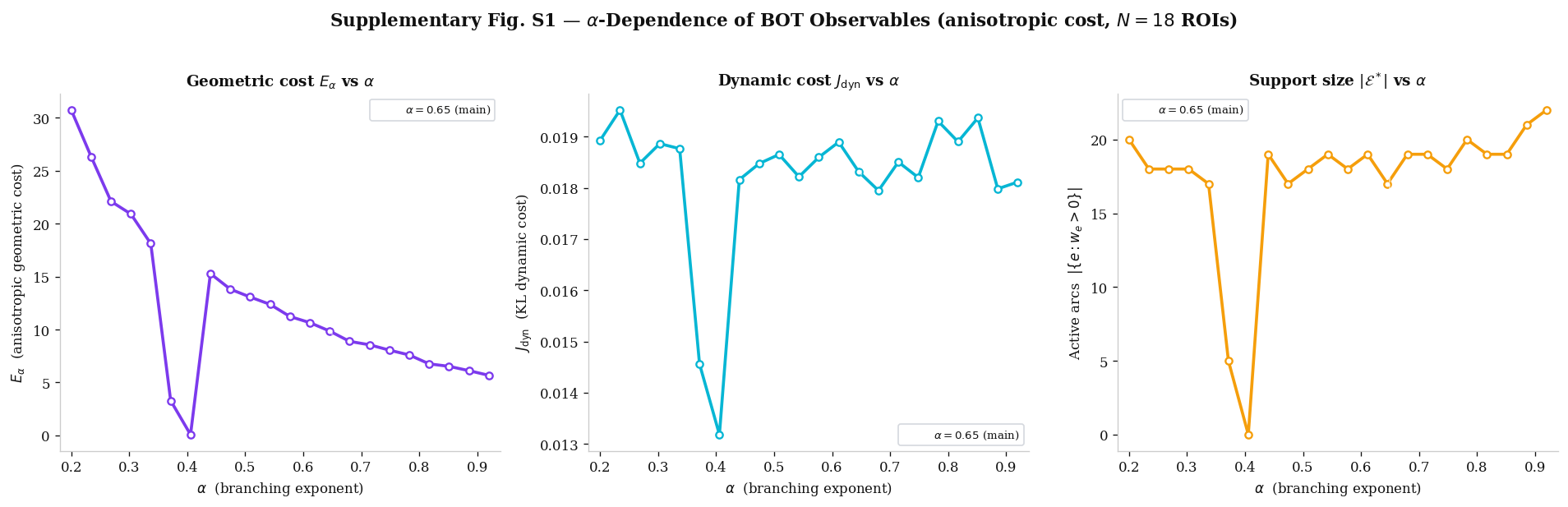}
\caption{\textbf{Supplementary Fig. S1. \(\alpha\)-dependence of anisotropic BOT observables.}
Left: anisotropic geometric cost \(E_\alpha\) versus branching exponent \(\alpha\). Middle: dynamic cost \(J_{\mathrm{dyn}}\) versus \(\alpha\). Right: support size \(|e^\ast|\) versus \(\alpha\). The vertical dashed line marks the main-text value \(\alpha=0.65\). The main stable regime is interrupted by a narrow near-degenerate transition around \(\alpha\approx 0.4\), where the support collapses and the observables show abrupt excursions.}
\label{fig:S1}
\end{figure*}

\subsection*{Representative anisotropic maps across branching regimes}

Supplementary Fig.~S2 shows anisotropic reaction maps for four representative values of \(\alpha\): \(0.25\), \(0.45\), \(0.65\), and \(0.85\). The figure makes the branching logic visually transparent. Lower \(\alpha\) values produce stronger aggregation and a smaller number of dominant high-load corridors, whereas larger \(\alpha\) values support more distributed routing across the graph. At the same time, the large-scale organisation of the backbone remains recognisable across the stable range, indicating that the main pathways are not completely rearranged by modest changes in branching strength.

This is important for interpretation. The branching exponent does not merely rescale flux magnitudes; it changes the architecture itself. The reaction map should therefore be understood as belonging to a parameterised family of economical routing backbones, within which the main-text solution at \(\alpha=0.65\) occupies an intermediate and structurally stable regime.

\newpage

\begin{figure*}[!t]
    \centering
    \includegraphics[width=\textwidth]{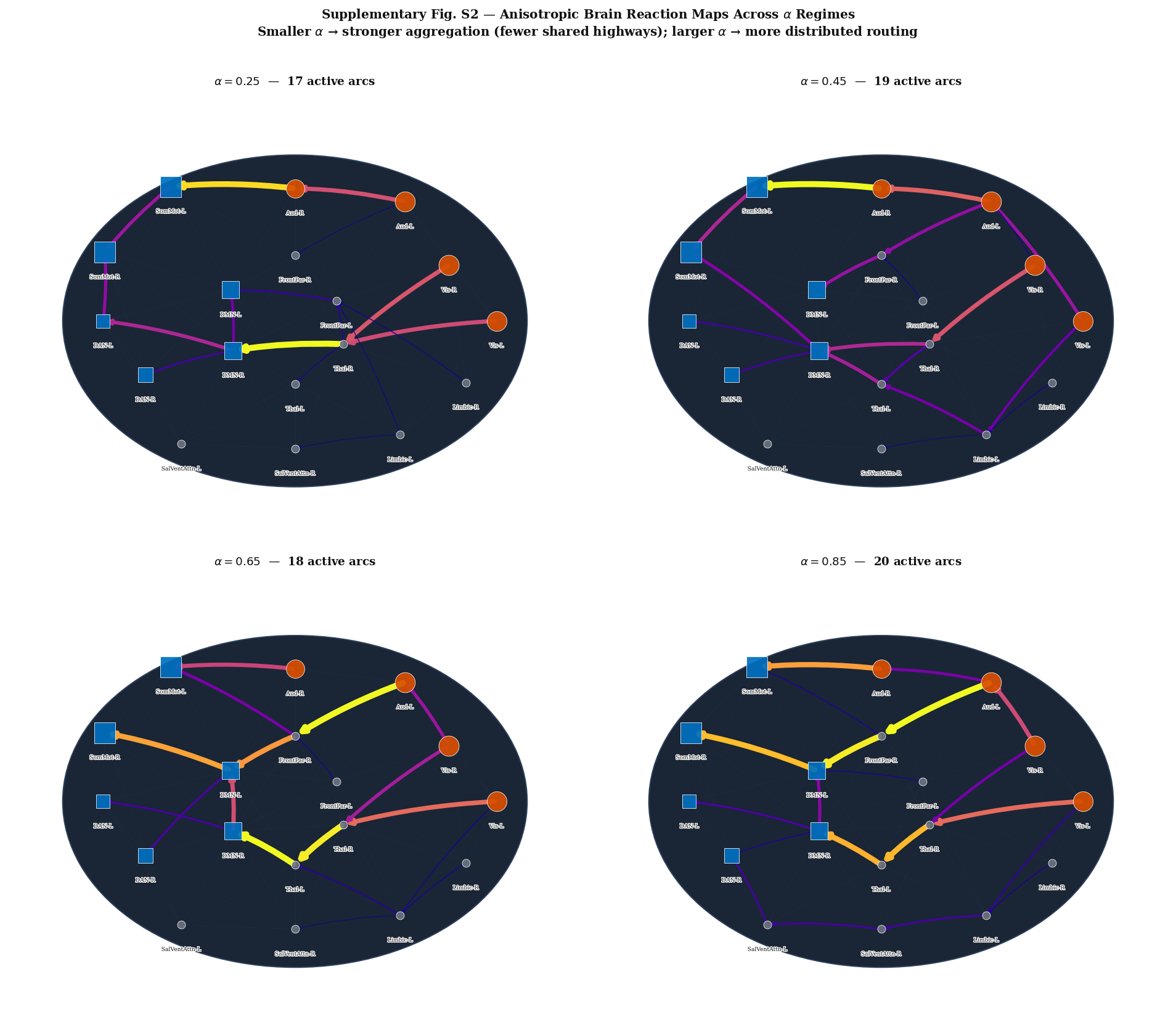}
\caption{\textbf{Supplementary Fig. S2. Anisotropic brain reaction maps across branching regimes.}
Reaction maps inferred at \(\alpha=0.25\), \(0.45\), \(0.65\), and \(0.85\). Smaller \(\alpha\) produces stronger aggregation and fewer shared highways, whereas larger \(\alpha\) supports more distributed routing. The backbone remains structurally interpretable across the stable parameter range.}
\label{fig:S2}
\end{figure*}

\begin{table}[H]
\centering
\caption{\textbf{Supplementary Table S1. Representative anisotropic supports across branching regimes.}}
\label{tab:S1}
\begin{tabular}{ccc}
\toprule
\(\alpha\) & Active arcs & Qualitative regime \\
\midrule
0.25 & 17 & strong aggregation \\
0.45 & 19 & post-transition stable regime \\
0.65 & 18 & main-text reference regime \\
0.85 & 20 & more distributed routing \\
\bottomrule
\end{tabular}
\end{table}

\section*{Supplementary Note 2. Comparison with non-branched baselines}

A critical question for the interpretation of the main results is whether the inferred shared corridors are a distinctive outcome of anisotropic branched transport or whether they would also emerge under simpler routing principles. Supplementary Fig.~S3 addresses this by comparing the anisotropic BOT solution with two non-branched surrogates defined on the same source--target pair and the same candidate graph: an anisotropic linear-cost transport surrogate (\(\alpha=1\)) and a greedy shortest-path surrogate.

\subsection*{What the baseline comparison shows}

The comparison reveals a clear structural hierarchy. The anisotropic BOT solution concentrates transport onto a relatively small set of shared relay corridors. The linear-cost surrogate weakens this concentration and spreads transport over a broader active support. The shortest-path surrogate fragments routing even further, favouring direct pairwise routes over mesoscale redistribution. Thus, the \emph{shared-highway architecture} that motivates the main manuscript is not simply a graph sparsity effect, nor a generic consequence of anisotropic weights. It is the specific consequence of combining anisotropy with a concave transport objective.

This distinction matters because it clarifies the conceptual gain over standard transport or routing baselines. The branched formulation does not only connect sources to sinks cheaply; it explains why separate signal contributions may become temporarily co-routed through common mesoscale corridors before diverging toward reaction-producing regions.

\begin{figure*}[!t]
    \centering
    \includegraphics[width=\textwidth]{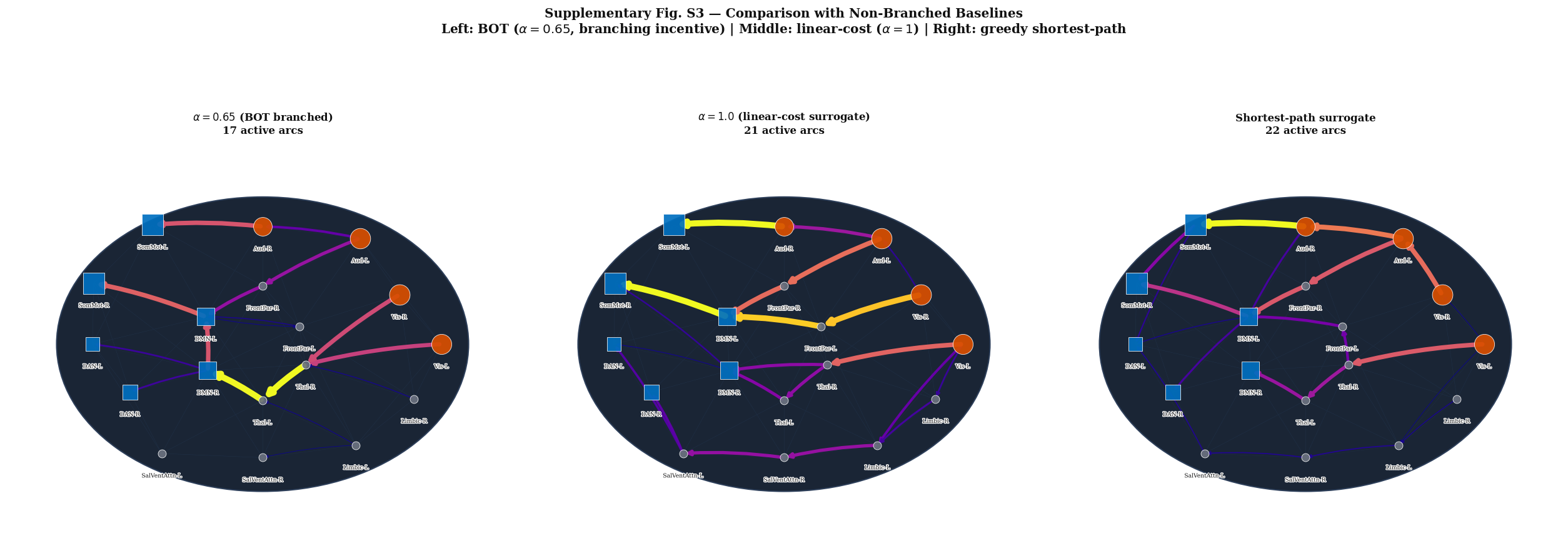}
\caption{\textbf{Supplementary Fig. S3. Comparison with non-branched baselines.}
Left: anisotropic BOT with \(\alpha=0.65\). Middle: anisotropic linear-cost surrogate (\(\alpha=1\)). Right: greedy shortest-path surrogate. The branched model concentrates transport onto a smaller number of relay corridors, whereas the non-branched surrogates distribute flux more diffusely and reduce the emergence of shared highways.}
\label{fig:S3}
\end{figure*}

\section*{Supplementary Note 3. Relay-region statistics and mesoscale interpretation}

The main manuscript argues that the anisotropic solution is biologically interpretable because the strongest relay regions are not arbitrary graph intermediates. Supplementary Fig.~S5 and Supplementary Table~S2 make this point quantitative.

\subsection*{Dominant relay systems}

The left panel of Supplementary Fig.~S5 ranks the top relay nodes by relay score. The strongest relays are the bilateral thalamic nodes, followed by bilateral default-mode and frontoparietal nodes. The right panel shows node-level incoming and outgoing flux, distinguishing relay, source, and sink nodes. Relay nodes lie near the diagonal because they both receive and redistribute substantial flow, whereas source- and sink-dominated nodes are displaced toward input- or output-heavy regimes.

This structure is highly informative. A direct-routing model might have concentrated almost all load on dominant source and target regions themselves. Instead, the branched anisotropic solution identifies intermediate systems that act as economical redistribution bottlenecks. In the present synthetic setting, the bilateral thalamic, default-mode, and frontoparietal nodes therefore emerge as the principal mesoscale bridge systems of the inferred reaction map.

\begin{figure*}[!t]
    \centering
    \includegraphics[width=\textwidth]{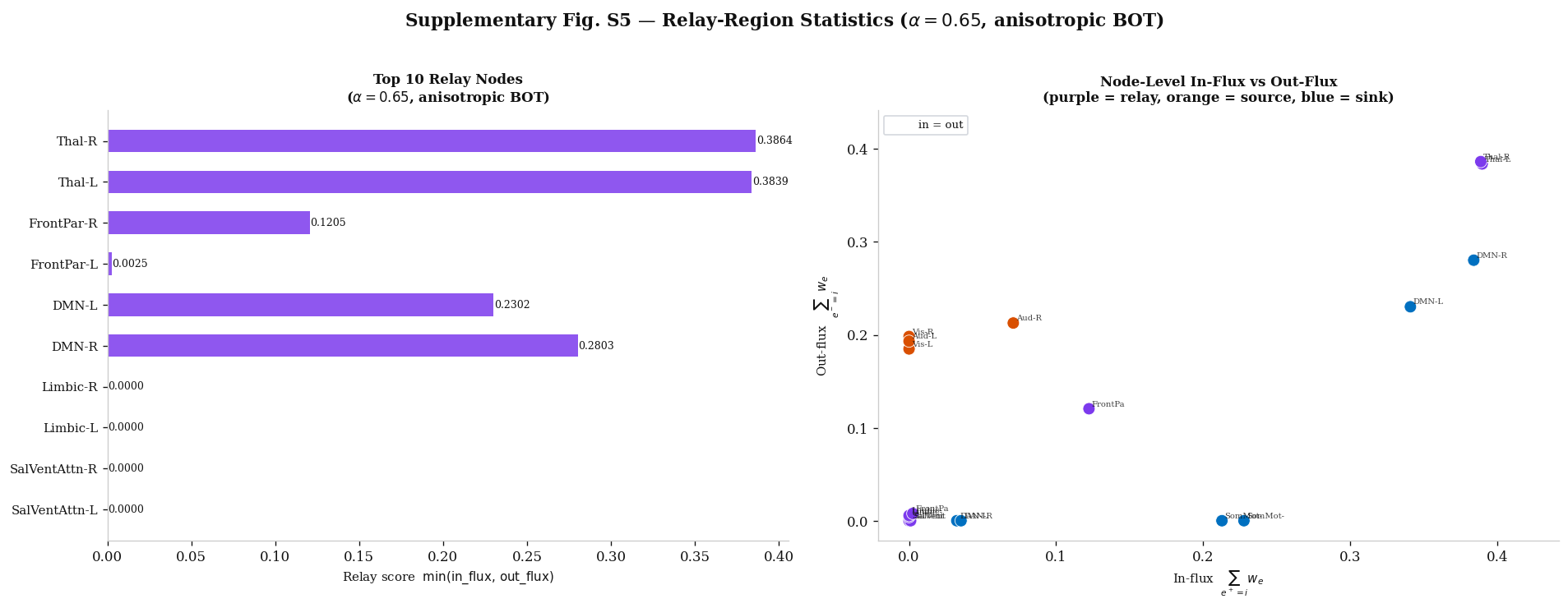}
\caption{\textbf{Supplementary Fig. S5. Relay-region statistics for the anisotropic solution at \(\alpha=0.65\).}
Left: ranking of the top relay nodes by relay score. Right: node-level incoming and outgoing flux, with relay nodes shown in purple, source-dominated nodes in orange, and sink-dominated nodes in blue. The dominant relay structure is led by bilateral thalamic, default-mode, and frontoparietal nodes.}
\label{fig:S5}
\end{figure*}

\begin{table}[H]
\centering
\caption{\textbf{Supplementary Table S2. Dominant relay nodes for the anisotropic solution at \(\alpha=0.65\).}}
\label{tab:S2}
\begin{tabular}{lc}
\toprule
Node & Relay score \\
\midrule
Thal-R & 0.3864 \\
Thal-L & 0.3839 \\
DMN-R & 0.2803 \\
DMN-L & 0.2302 \\
FrontPar-R & 0.1205 \\
FrontPar-L & 0.0025 \\
Limbic-L & 0.0000 \\
Limbic-R & 0.0000 \\
SalVentAttn-R & 0.0000 \\
SalVentAttn-L & 0.0000 \\
\bottomrule
\end{tabular}
\end{table}

The dominance of bilateral thalamic and default-mode nodes is especially relevant because it indicates that the inferred architecture is not merely a shortest transfer from sensory entry to motor output. Instead, the model selects a hierarchical backbone in which integrative bridge systems mediate aggregation and redistribution at a mesoscale level.

\section*{Supplementary Note 4. Proof-of-principle scalability analysis}

One of the most natural concerns raised by the main manuscript is computational scalability: can the framework, demonstrated here on \(18\) regions of interest, plausibly extend to larger connectomic resolutions? Supplementary Fig.~S4 provides a proof-of-principle answer.

\subsection*{Observed scaling over the tested graph family}

The left panel reports solver runtime as a function of graph size, while the right panel reports the growth in the number of directed arcs under the adopted \(k\)-nearest-neighbour graph construction. Over the tested synthetic graph family, the displayed log--log fits suggest a runtime scaling approximately proportional to \(N^{2.02}\) and an edge-count scaling approximately proportional to \(N^{1.09}\).

These fits should not be overinterpreted as universal asymptotic laws. Their role is different: they show that over the explored regime the candidate graph remains sparse and the solver displays roughly near-quadratic growth. For the present proof-of-principle article, this is enough to support the claim that the framework is not inherently restricted to the \(18\)-ROI demonstration.

\subsection*{Interpretation for larger connectomes}

A connectome with \(200+\) nodes would still be substantially more demanding, especially once repeated solves across \(\alpha\) and repeated dynamic-cost evaluations are taken into account. However, the current scaling curve supports a reasonable path forward. The most natural routes to improved scalability are:
\begin{enumerate}[label=(\alph*)]
    \item anatomically informed sparsification of candidate arcs,
    \item warm-start or continuation strategies across neighbouring \(\alpha\) values,
    \item parallel evaluation of graph-dependent dynamic costs,
    \item alternative large-scale solvers tailored to sparse nonlinear transport problems.
\end{enumerate}

Thus, Supplementary Fig.~S4 should be read as evidence that the current implementation is a \emph{transparent variational proof-of-principle} with plausible computational extensions, rather than as a claim that the current SLSQP-based code is already the final large-connectome production solver.

\begin{figure*}[!t]
    \centering
    \includegraphics[width=\textwidth]{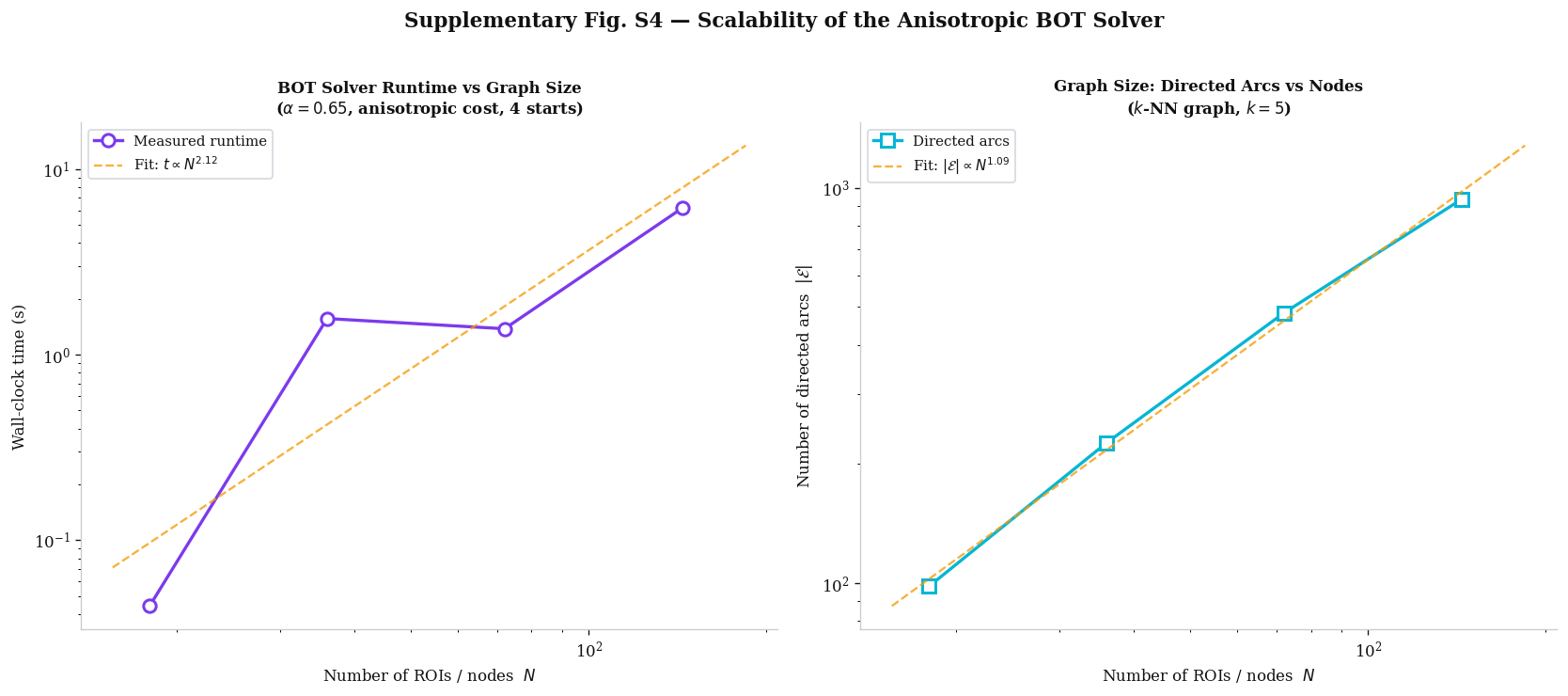}
\caption{\textbf{Supplementary Fig. S4. Scalability of the anisotropic BOT solver.}
Left: measured BOT runtime versus graph size, with fitted power-law scaling approximately proportional to \(N^{2.02}\). Right: number of directed arcs versus number of nodes, with fitted scaling approximately proportional to \(N^{1.09}\). The graph family remains sparse, while runtime grows near-quadratically over the tested range.}
\label{fig:S4}
\end{figure*}

\section*{Supplementary Note 5. Strengthened interpretation of the supplementary analyses}

Taken together, the supplementary analyses reinforce the main manuscript in four decisive ways.

First, the \(\alpha\)-sensitivity analysis shows that the anisotropic reaction map belongs to a structured family of routing architectures, rather than depending on one isolated parameter choice. Second, the baseline comparison makes explicit that the shared relay corridors of the main manuscript are a genuine signature of branched transport and are not recovered by standard non-branched alternatives. Third, the relay statistics demonstrate that the strongest intermediates are biologically interpretable mesoscale systems rather than anonymous connectors. Fourth, the scalability analysis clarifies that the present work is best understood as a variational proof-of-principle that opens the way to larger connectomic implementations.

The purpose of this Supplementary Information is therefore not merely additive. It sharpens the manuscript’s scientific position. The main result is not just that one can solve a synthetic transport problem on a small graph, but that a multimodal anisotropic branched-transport framework can reveal interpretable routing backbones, distinguish them from simpler baselines, and organise them into a geometric--dynamic trade-off landscape whose structure persists across a stable branching regime.

\end{document}